\documentclass[english,11pt]{article}
\usepackage[latin1]{inputenc}
\usepackage{amsmath,amsthm,amssymb,babel}

\newtheorem{teo}{Theorem}

\newtheorem{lemma}{Lemma}

\def\proof{{\it Proof.}\ }
\def\endproof{\hfill $\Box$\par\vskip3mm}

\def\eq#1{(\ref{#1})}

\def\neweq#1{\begin{equation}\label{#1}}
\def\endeq{\end{equation}}

\def\phi{\varphi}
\def\RR{{\mathbb R} }
\def\NN{{\mathbb N} }

\def\di{\displaystyle}
\def\ri{\rightarrow}

\date{}

\title{\sc A continuous spectrum for nonhomogeneous differential operators in
 Orlicz-Sobolev spaces\thanks{
Correspondence address: Vicen\c{t}iu R\u{a}dulescu, Department of
Mathematics, University of Craiova,  200585 Craiova, Romania. E-mail:
{\tt
vicentiu.radulescu@math.cnrs.fr}}}
\author{\sc Mihai Mih\u ailescu$\,^{a,b}$ and Vicen\c{t}iu R\u{a}dulescu$\,^{a,c}$\\
\small
$^a\,$Department of Mathematics, University of Craiova,  200585 Craiova,
Romania\\
\small $^b\,$Department of Mathematics, Central European University, 1051 Budapest, Hungary\\
\small $^c\,$Institute of Mathematics ``Simion Stoilow" of the Romanian Academy,\\
\small P.O. Box 1-764, 014700 Bucharest, Romania\\
\small
E-mail addresses: {\tt mmihailes@yahoo.com}\qquad {\tt
vicentiu.radulescu@math.cnrs.fr}}

\begin{document}
\baselineskip16pt \maketitle \noindent{\small{\sc Abstract}. We
study the nonlinear eigenvalue problem $-{\rm div}(a(|\nabla
u|)\nabla u)=\lambda|u|^{q(x)-2}u$ in $\Omega$, $u=0$ on
$\partial\Omega$, where $\Omega$ is a bounded open set in $\RR^N$
with smooth boundary, $q$ is a continuous function, and $a$ is a nonhomogeneous potential.
We establish sufficient conditions on $a$ and $q$ such that
 the above nonhomogeneous quasilinear problem has continuous families of
eigenvalues. The proofs rely   on elementary variational arguments. The abstract results of this paper
are illustrated by the cases $a(t)=t^{p-2}\log (1+t^r)$ and $a(t)= t^{p-2} [\log (1+t)]^{-1}$.
}

\section{Introduction and preliminary results}
Let $\Omega$ be a bounded domain in $\RR^N$ ($N\geq 3$) with smooth boundary
$\partial\Omega$.
In this paper we are concerned with the following eigenvalue problem:
\begin{equation}\label{1}
\left\{\begin{array}{lll}
-{\rm div}(a(|\nabla u|)\nabla u)=\lambda|u|^{q(x)-2}u, &\mbox{for}&
x\in\Omega\\
u=0, &\mbox{for}& x\in\partial\Omega\,.
\end{array}\right.
\end{equation}
We assume that the function $a:(0,\infty)\ri\RR$ is such that the mapping $\phi:\RR\rightarrow\RR$ defined by
$$\phi(t)=\left\{\begin{array}{lll}
a(|t|)t, &\mbox{for}&
t\neq 0\\
0, &\mbox{for}& t=0\,,
\end{array}\right.$$
is an odd, increasing homeomorphism from $\RR$ onto $\RR$. We also suppose throughout this paper that
 $\lambda>0$ and $q:\overline\Omega\ri (0,\infty )$ is a continuous function.

Since the
operator in the divergence form is nonhomogeneous we
introduce an Orlicz-Sobolev space setting for problems of this type.
On the other hand, the term arising in the right hand side
of equation \eq{1} is also nonhomogeneous and its particular form
appeals to a suitable variable exponent Lebesgue space setting.

We point out that eigenvalue problems involving quasilinear nonhomogeneous problems in Orlicz-Sobolev spaces
were studied in \cite{Gar} but in a different framework. In what concerns the case when $a(|\nabla u|)=
|\nabla u|^{q(x)-2}$, problem \eq{1} was studied by
Fan {\it et al.} in \cite{fann,FZZ} who established the existence of
a sequence of eigenvalues, by means of the Ljusternik-Schnirelmann critical point theory.  Denoting by $\Lambda$ the 
set of all
nonnegative eigenvalues, Fan, Zhang and Zhao showed that $\sup\Lambda=+\infty$ and they pointed out that
only under additional assumptions we have $\inf\Lambda>0$. We remark that in the homogeneous case
corresponding the $p$-Laplace
operator ( $p(x)\equiv p$) we always have $\inf\Lambda>0$. A different approach of the eigenvalue problem
\eq{1} corresponding to $a(|\nabla u|)=
|\nabla u|^{p(x)-2}$ and $p(x)\not=q(x)$ is given in Mih\u ailescu and R\u adulescu \cite{pams}.

We first recall some basic facts about Orlicz spaces.
Define
$$\Phi(t)=\int_0^t\phi(s)\;ds,\;\;\;\Phi^\star(t)=\int_0^t\phi^{-1}(s)\;ds,\qquad \mbox{for all}\ t\in\RR.$$
We observe that $\Phi$ is a {\it Young function}, that is, $\Phi (0)=0$, $\Phi$ is convex, and $\lim_{x\ri\infty}\Phi 
(x)=+\infty$.
Furthermore, since $\Phi (x)=0$ if and only if $x=0$, $\lim_{x\ri 0}\Phi (x)/x=0$, and $\lim_{x\ri \infty}\Phi 
(x)/x=+\infty$, then
$\Phi$ is called an {\it $N$--function}.
The function $\Phi^\star$ is called the {\it complementary} function of $\Phi$ and it satisfies
$$\Phi^\star (t)=\sup\{st-\Phi (s);\ s\geq 0\},\qquad\mbox{for all $t\geq 0$}\,.$$
 We observe that $\Phi^\star$ is also an $N$--function and the following Young's inequality holds true:
 $$st\leq\Phi (s)+\Phi^\star (t),\qquad\mbox{for all $s,t\geq 0$}\,.$$

The Orlicz space $L_\Phi(\Omega)$ defined by the $N$--function $\Phi$
(see \cite{AHed,A,Clem1}) is the space of measurable functions $u:\Omega\ri\RR$ such that
$$\|u\|_{L_\Phi}:=\sup\left\{\int_\Omega uv\;dx;\ \int_\Omega\Phi^\star (|v|)\;dx\leq 1\right\}<\infty\,.$$
Then $(L_\Phi(\Omega),\|\,\cdot\, \|_{L_\Phi} )$ is a Banach space whose norm is equivalent to the
 Luxemburg norm
$$\|u\|_\Phi :=\inf\left\{k>0;\ \int_\Omega\Phi\left(\frac{u(x)}{k}
\right)\;dx\leq 1\right\}.$$
For Orlicz spaces the H\"older's inequality reads as follows (see \cite[Inequality 4, p. 79]{rao}):
$$\int_\Omega uvdx\leq 2\,\|u\|_{L_\Phi}\, \|v\|_{L_{\Phi^\star}}\qquad\mbox{for all $u\in L_\Phi(\Omega)$
and $v\in L_{\Phi^\star}(\Omega)$}\,.$$

We denote by $W_0^1L_\Phi(\Omega)$ the corresponding Orlicz-Sobolev space for problem \eq{1},
equipped with the norm
$$\|u\|=\||\nabla u|\|_\Phi$$
(see \cite{A,Clem1,G}). The space $W_0^1L_\Phi(\Omega)$ is also a Banach space.

In this paper we assume that $\Phi$ and $\Phi^\star$ satisfy the $\Delta_2$-condition (at infinity), namely
$$1<\liminf\limits_{t\rightarrow\infty}\frac{t\phi(t)}{\Phi(t)}
\leq\limsup\limits_{t>0}\frac{t\phi(t)}{\Phi(t)}<\infty.$$
Then $L_\Phi(\Omega)$ and $W_0^1L_\Phi(\Omega)$ are reflexive Banach spaces.

Now we introduce the Orlicz-Sobolev conjugate $\Phi_\star$ of $\Phi$, defined as
$$\Phi_\star^{-1}(t)=\int_0^t\frac{\Phi^{-1}(s)}
{s^{(N+1)/N}}\;ds.$$
We assume that
\begin{equation}\label{2}
\lim_{t\rightarrow 0}\int_t^1\frac{\Phi^{-1}(s)}{s^{(N+1)/N}}\;ds<\infty,\;\;\;{\rm and}\;\;\;
\lim_{t\rightarrow\infty}\int_1^t\frac{\Phi^{-1}(s)}{s^{(N+1)/N}}\;ds=\infty.
\end{equation}
Finally, we define
$$p_0:=\inf_{t>0}\frac{t\phi(t)}{\Phi(t)}\;\;{\rm and}\;\;p^0:=\sup_{t>0}\frac{t\phi(t)}{\Phi(t)}.$$

Next, we recall some background facts concerning the variable exponent Lebesgue spaces. For more details we refer to 
the book by Musielak
\cite{M} and the papers by Acerbi {\it at al.}, Edmunds {\it et al.} \cite{edm, edm2, edm3}, Kovacik and  
R\'akosn\'{\i}k \cite{KR},
Mih\u ailescu and R\u adulescu \cite{RoyalSoc}, Samko and Vakulov \cite{samko}, Zhikov \cite{zhikov}.

Set
$$C_+(\overline\Omega)=\{h;\;h\in C(\overline\Omega),\;h(x)>1\;{\rm
for}\;
{\rm all}\;x\in\overline\Omega\}.$$
For any $h\in C_+(\overline\Omega)$ we define
$$h^+=\sup_{x\in\Omega}h(x)\qquad\mbox{and}\qquad h^-=
\inf_{x\in\Omega}h(x).$$
For any $q(x)\in C_+(\overline\Omega)$ we define the variable exponent
Lebesgue space $L^{q(x)}(\Omega)$ (see \cite{KR}).
On $L^{q(x)}(\Omega)$ we define the Luxemburg norm by the formula
$$|u|_{q(x)}=\inf\left\{\mu>0;\;\int_\Omega\left|
\frac{u(x)}{\mu}\right|^{q(x)}\;dx\leq 1\right\}.$$ We remember that
the variable exponent Lebesgue spaces are separable and reflexive
Banach spaces. If $0 <|\Omega|<\infty$ and $q_1$, $q_2$ are variable
exponents so that $q_1(x) \leq q_2(x)$ almost everywhere in $\Omega$
then there exists the continuous embedding
$L^{q_2(x)}(\Omega)\hookrightarrow L^{q_1(x)}(\Omega)$.

If $(u_n)$, $u\in L^{q(x)}(\Omega)$ then the following relations
hold true
\begin{equation}\label{3}
|u|_{q(x)}>1\;\;\;\Rightarrow\;\;\;|u|_{q(x)}^{q^-}\leq\int_\Omega|u|^{q(x)}\;dx
\leq|u|_{q(x)}^{q^+}
\end{equation}
\begin{equation}\label{4}
|u|_{q(x)}<1\;\;\;\Rightarrow\;\;\;|u|_{q(x)}^{q^+}\leq
\int_\Omega|u|^{q(x)}\;dx\leq|u|_{q(x)}^{q^-}
\end{equation}
\begin{equation}\label{4biss}
|u_n-u|_{q(x)}\rightarrow 0\;\;\;\Leftrightarrow\;\;\;\int_\Omega|u_n-u|^{q(x)}\;dx
\rightarrow 0.
\end{equation}

\section{The main results}
We say that $\lambda\in\RR$ is an eigenvalue of problem \eq{1} if
there exists $u\in W_0^1L_\Phi(\Omega)
\setminus\{0\}$ such that
$$\int_{\Omega}a(|\nabla u|)\nabla u\nabla v\;dx-\lambda\int_{\Omega}
|u|^{q(x)-2}uv\;dx=0,$$ for all $v\in W_0^1L_\Phi(\Omega)$. We point
out that if $\lambda$ is an eigenvalue of  problem \eq{1} then the
corresponding $u\in W_0^1L_\Phi(\Omega)\setminus\{0\}$ is a  weak
solution of \eq{1}, called an eigenvector of equation \eq{1}
corresponding to the eigenvalue $\lambda$.

Our first main result shows that, in certain circumstances, any
positive and sufficiently small $\lambda$ is an eigenvalue of \eq{1}.

\begin{teo}\label{t1}
Assume that relation \eq{2} is fulfilled and furthermore
\begin{equation}\label{6}
1<\inf_{x\in\Omega}q(x)<p_0\,,
\end{equation}
and
\begin{equation}\label{6biss}
\lim_{t\rightarrow\infty}\di\frac{|t|^{q^+}}
{\Phi_\star(kt)}=0,\;\;\; {\rm for}\; {\rm all}\; k>0.
\end{equation}
Then there exists $\lambda^\star>0$
such that any $\lambda\in(0,\lambda^\star)$ is an eigenvalue for problem \eq{1}.
\end{teo}

The above result implies
$$\inf_{u\in W_0^1L_\Phi(\Omega)\setminus\{0\}}\frac{\di\int_\Omega
\Phi(|\nabla u|)\;dx}{\di\int_\Omega|u|^{q(x)}\;dx}=0.$$

The second main result of this paper asserts that in certain cases the set of eigenvalues may coincide with the {\it 
whole} positive semiaxis.

\begin{teo}\label{t2}
Assume that relations \eq{2} and \eq{6biss} are fulfilled and furthermore
\begin{equation}\label{stea2}
\sup_{x\in\Omega}q(x)<p_0.
\end{equation}
Then every $\lambda>0$ is an eigenvalue for problem \eq{1}. Moreover, for any $\lambda>0$
there exists a sequence of eigenvectors $\{u_n\}\subset E$ such that $\lim_{n\rightarrow\infty}u_n=0$
in $W_0^1L_\Phi(\Omega)$.
\end{teo}

\noindent{\bf Remark 1.} Relations \eq{2} and \eq{6biss} enable us
to apply Theorem 2.2 in \cite{Gar} (see also Theorem 8.33 in
\cite{A}) in order to obtain that $W_0^1L_\Phi(\Omega)$ is compactly
embedded in $L^{q^+}(\Omega)$. This fact combined with the
continuous embedding of $L^{q^+}(\Omega)$ in $L^{q(x)}(\Omega)$
ensures that $W_0^1L_\Phi(\Omega)$ is compactly embedded in
$L^{q(x)}(\Omega)$.
\smallskip

\noindent{\bf Remark 2.} The conclusion of Theorems \ref{t1} and
\ref{t2} still remains valid if we replace the hypothesis \eq{6biss}
in Theorems \ref{t1} and \ref{t2} by the following relation
\begin{equation}\label{66}
N<p_0<\liminf_{t\rightarrow\infty}\frac{\log(\Phi(t))}{\log(t)}.
\end{equation}
Indeed, using Lemma D.2 in \cite{Clem2}, it follows that $W_0^1L_\Phi(\Omega)$ is continuously embedded in
$W_0^{1,p_0}(\Omega)$. On the other hand, since we assume $p_0>N$, we deduce that $W_0^{1,p_0}(\Omega)$ is compactly 
embedded
in $C(\overline\Omega)$. Thus, we obtain
 that $W_0^1L_\Phi(\Omega)$ is compactly embedded in $C(\overline\Omega)$. Since
$\Omega$ is bounded it follows that $W_0^1L_\Phi(\Omega)$ is continuously embedded in $L^{q(x)}(\Omega)$.

\section{Proof of Theorem \ref{t1}}
Let $E$ denote the Orlicz-Sobolev space $W_0^1L_\Phi(\Omega)$.

For any $\lambda>0$ the energy functional $J_\lambda:E\rightarrow\RR$ corresponding to problem \eq{1} is defined by
$$J_\lambda(u)=\int_{\Omega}\Phi(|\nabla u|)\;dx-\lambda
\int_\Omega\frac{1}{q(x)}|u|^{q(x)}\;dx.$$
Standard arguments imply that
$J_\lambda\in C^1(E,\RR)$ and
$$\langle J_\lambda^{'}(u),v\rangle=\int_{\Omega}a(|\nabla u|)
\nabla u\nabla v\;dx-\lambda\int_\Omega|u|^{q(x)-2}uv\;dx,$$ for all
$u,\;v\in E$. Thus the weak solutions of \eq{1} coincide with the
critical points of $J_\lambda$. If such a weak solution exists and
is nontrivial then the corresponding $\lambda$ is an eigenvalue of
problem \eq{1}.

\begin{lemma}\label{l1}
There is some $\lambda^\star>0$ such that for any $\lambda\in(0,\lambda^\star)$
there exist $\rho$, $\alpha>0$  such that $J_\lambda(u)\geq\alpha>0$
for any $u\in E$ with $\|u\|=\rho$.
\end{lemma}

\proof By the definition of $p^0$ and since $\frac{{\mathrm d}}{{\mathrm d}\tau}\left(\tau^{p_0}\Phi (t/\tau)\right)
\geq 0$ we obtain
$$\Phi(t)\geq\tau^{p^0}\Phi(t/\tau),\;\;\;\forall\; t>0\;{\rm and}\;
\tau\in(0,1]\,,$$
(see page 44 in \cite{Clem1}). Combining this fact with Proposition 6 in \cite[page 77]{rao}
we find that
\begin{equation}\label{7}
\int_\Omega\Phi(|\nabla u(x)|)\;dx\geq\|u\|^{p^0},\;\;\;\forall\;
u\in E\;{\rm with}\;\|u\|<1.
\end{equation}

On the other hand, since $E$ is continuously embedded in $L^{q(x)}(\Omega)$,
there exists a positive constant $c_1$ such that
\begin{equation}\label{8}
|u|_{q(x)}\leq c_1\|u\|,\;\;\;\forall\;u\in E.
\end{equation}
We fix $\rho\in(0,1)$ such that $\rho<1/c_1$. Then  relation  \eq{8} implies
\begin{equation}\label{9}
|u|_{q(x)}<1,\;\;\;\forall\;u\in E,\;{\rm with}\;\|u\|=\rho.
\end{equation}
Furthermore, relation \eq{4} yields
\begin{equation}\label{10}
\int_\Omega|u|^{q(x)}\;dx\leq|u|_{q(x)}^{q^-},\;\;\;\forall\;u\in E,\;{\rm with}\;\|u\|=\rho.
\end{equation}
Relations \eq{8} and \eq{10} imply
\begin{equation}\label{10biss}
\int_\Omega|u|^{q(x)}\;dx\leq c_1^{q^-}\|u\|^{q^-},\;\;\;\forall\;u\in E,\;{\rm with}\;\|u\|=\rho.
\end{equation}
Taking into account  relations \eq{7}, \eq{4} and \eq{10biss} we deduce that
for any $u\in E$ with $\|u\|=\rho$ the following inequalities hold true
$$
J_\lambda(u)\geq\|u\|^{p^0}-\frac{\lambda}{q^-}
\int_\Omega|u|^{q(x)}\;dx
=\rho^{q^-}\left(\rho^{p^0-q^-}-\frac{\lambda}{q^-}c_1^{q^-}\right).
$$
We point out that by relation \eq{6} and the definition of $p^0$ we have
$q^-<l\leq p^0$.
By the above inequality we remark that if we define
\begin{equation}\label{11}
\lambda^\star=\frac{\rho^{p^0-q^-}}{2}\cdot\frac{q^-}{c_1^{q^-}}
\end{equation}
then for any $\lambda\in(0,\lambda^\star)$ and any $u\in E$ with $\|u\|=\rho$ there
exists $\alpha=\frac{\rho^{p^0}}{2}>0$ such that
$$J_\lambda(u)\geq\alpha>0.$$
The proof of Lemma \ref{l1} is complete.  \endproof

\begin{lemma}\label{l2}
There exists $\phi\in E$ such that $\phi\geq 0$,
$\varphi\neq 0$ and
$J_\lambda(t\phi)<0$,
for $t>0$ small enough.
\end{lemma}

\proof
Assumption \eq{6} implies that $q^-<p_0$. Let $\epsilon_0>0$ be such that $q^-+\epsilon_0<p_0$.
On the other hand, since $q\in C(\overline\Omega)$ it follows that there exists an open set
$\Omega_0\subset\Omega$ such that $|q(x)-q^-|<\epsilon_0$ for all $x\in\Omega_0$. Thus, we
conclude that $q(x)\leq q^-+\epsilon_0<p_0$ for all $x\in\Omega_0$.

Let $\phi\in C_0^\infty(\Omega)$ be such that ${\rm supp}(\phi)\supset\overline\Omega_0$,
$\phi(x)=1$ for all $x\in\overline\Omega_0$ and $0\leq\phi\leq 1$ in $\Omega$.

We also point out that there exists $t_0\in(0,1)$ such that for any $t\in(0,t_0)$ we have
$$\|t|\nabla\phi|\|_\Phi=t\|\phi\|<1.$$
Taking into account all the above information and using Lemma C.9 in \cite{Clem2} we have
\begin{eqnarray*}
J_\lambda(t\phi)&=&\int_{\Omega}\Phi(t|\nabla\phi(x)|)\;dx-\lambda
\int_\Omega\frac{t^{q(x)}}{q(x)}|\phi|^{q(x)}\;dx\\
&\leq&\int_{\Omega}\Phi(t|\nabla\phi(x)|)\;dx-\frac{\lambda}{q^+}
\int_\Omega t^{q(x)}|\phi|^{q(x)}\;dx\\
&\leq&\int_{\Omega}\Phi(t|\nabla\phi(x)|)\;dx-\frac{\lambda}{q^+}
\int_{\Omega_0} t^{q(x)}|\phi|^{q(x)}\;dx\\
&\leq&t^{p_0}\|\phi\|^{p_0}-\frac{\lambda\cdot
t^{q^-+\epsilon_0}}{q^+}|\Omega_0|,
\end{eqnarray*}
for any $t\in(0,1)$, where $|\Omega_0|$  denotes the Lebesgue measure of $\Omega_0$.
Therefore
$$J_\lambda(t\phi)<0$$
for $t<\delta^{1/(p_0-q^--\epsilon_0)}$, where
$$0<\delta<\min\left\{t_0,\frac{\frac{\lambda}{q^+}|\Omega_0|}{\|\phi\|^{p_0}}\right\}.$$
The proof of Lemma \ref{l2} is complete.   \endproof

\medskip
{\sc Proof of Theorem \ref{t1}.}
Let $\lambda^\star>0$ be defined as in \eq{11} and $\lambda\in(0,\lambda^\star)$.
By Lemma \ref{l1} it follows that on the boundary of the  ball centered at the origin and of
radius $\rho$ in $E$, denoted by $B_\rho(0)$, we have
\begin{equation}\label{12}
\inf\limits_{\partial B_\rho(0)}J_\lambda>0.
\end{equation}
On the other hand, by Lemma \ref{l2}, there exists $\phi\in E$ such that
$J_\lambda(t\phi)<0$ for all $t>0$ small enough. Moreover, relations \eq{7}, \eq{10biss} and
\eq{4} imply that for any $u\in B_\rho(0)$ we have
$$J_\lambda(u)\geq\|u\|^{p^0}-\frac{\lambda}{q^-}c_1^{q^-}\|u\|^{q^-}\,.$$
It follows that
$$-\infty<\underline{c}:=\inf\limits_{\overline{B_\rho(0)}}J_\lambda<0.$$
We let now $0<\epsilon<\inf_{\partial B_\rho(0)}J_\lambda-\inf_{B_\rho(0)}J_\lambda$.
Applying Ekeland's variational principle \cite{E} to the functional
$J_\lambda:\overline{B_\rho(0)}\rightarrow\RR$,  we find
$u_\epsilon\in\overline{B_\rho(0)}$ such that
\begin{eqnarray*}
J_\lambda(u_\epsilon)&<&\inf\limits_{\overline{B_\rho(0)}}J_\lambda+\epsilon\\
J_\lambda(u_\epsilon)&<&J_\lambda(u)+\epsilon\cdot\|u-u_\epsilon\|,\;\;\;u\neq u_\epsilon.
\end{eqnarray*}
Since
$$J_\lambda(u_\epsilon)\leq\inf\limits_{\overline{B_\rho(0)}}J_\lambda+\epsilon\leq
\inf\limits_{B_\rho(0)} J_\lambda+\epsilon<\inf\limits_{\partial B_\rho(0)}J_\lambda\,,$$
we deduce that $u_\epsilon\in B_\rho(0)$. Now, we define $I_\lambda:
\overline{B_\rho(0)}\rightarrow\RR$ by $I_\lambda(u)=J_\lambda(u)+\epsilon\cdot
\|u-u_\epsilon\|$. It is clear that $u_\epsilon$ is a minimum point
of $I_\lambda$ and thus
$$\di\frac{I_\lambda(u_\epsilon+t\cdot v)-{I_\lambda}(u_\epsilon)}
{t}\geq 0$$
for small $t>0$ and any $v\in B_1(0)$. The
above relation yields
$$\di\frac{J_\lambda(u_\epsilon+t\cdot v)-J_\lambda(u_\epsilon)}{t}+
\epsilon\cdot\|v\|\geq 0.$$
Letting $t\rightarrow 0$ it follows that $\langle J_\lambda^{'}
(u_\epsilon),v\rangle+\epsilon\cdot\|v\|>0$ and we infer that
$\|J_\lambda^{'}(u_\epsilon)\|\leq\epsilon$.

We deduce that there exists a sequence
$\{w_n\}\subset B_\rho(0)$ such that
\begin{equation}\label{13}
J_\lambda(w_n)\rightarrow{\underline c}\;\;\;{\rm and}\;\;\;
J_\lambda^{'}(w_n)\rightarrow 0.
\end{equation}
It is clear that $\{w_n\}$ is bounded in $E$. Thus, there exists $w\in E$ such that, up  to a
subsequence, $\{w_n\}$ converges weakly to $w$ in $E$.
By Remark 2 we deduce that $E$ is compactly embeddded in
$L^{q(x)}(\Omega)$, hence $\{w_n\}$ converges strongly to $w$ in $L^{q(x)}(\Omega)$.
So, by relations \eq{4biss} and H\"{o}lder's inequality for variable exponent spaces (see e.g. \cite{KR}),
$$\lim\limits_{n\rightarrow\infty}\int_\Omega|w_n|^{q(x)}\;dx=\int_\Omega|w|^{q(x)}\;dx\;\;\;
{\rm and}\;\;\;\lim\limits_{n\rightarrow\infty}\int_\Omega|w_n|^{q(x)-2}w_nv\;dx=
\int_\Omega|w|^{q(x)-2}wv\;dx$$
for any $v\in E$.

We conclude that $w$ is a nontrivial weak solution for problem \eq{1} and thus any
$\lambda\in(0,\lambda^\star)$ is an eigenvalue of problem \eq{1}. Similar arguments as those
used on page 50 in \cite{Clem1} imply that $\{w_n\}$ converges strongly to $w$ in $E$. So, by \eq{13},
\begin{equation}\label{LLLLet11}
J_\lambda(w)=\underline c<0\;\;\;{\rm and}\;\;\;J_\lambda^{'}(w)=0.
\end{equation}

The proof of Theorem \ref{t1} is complete.\endproof

\section{Proof of Theorem \ref{t2}}
We still denote by $E$ the Orlicz-Sobolev space $W_0^1L_\Phi(\Omega)$. For any $\lambda>0$ let $J_\lambda$ be defined
as in the above section of the paper.

In order to prove Theorem \ref{t2} we apply to the functional $J_\lambda$ a symmetric version of
the mountain pass lemma, recently developed by Kajikia in \cite{K}. Before presenting the result in
\cite{K} we remember the following definition.
\smallskip

\noindent{\bf Definition 1.} Let $X$ be a real Banach space.
We say that a subset $A$ of $X$ is {\it symmetric} if
$u\in A$ implies $-u\in A$. For a closed symmetric set $A$ which
does not contain the origin, we define the {\it genus} $\gamma(A)$ of
$A$ as the smallest integer $k$ such that there exists an odd
continuous mapping from $A$ to $\RR^k\setminus\{0\}$. If there does
not exist such an integer $k$, we define $\gamma(A)=+\infty$. Moreover, we set
$\gamma(\emptyset)=0$. Finally, we denote by $\Gamma_k$ the family
$$\Gamma_k=\{A\subset X;\; 0\not\in A\; {\rm and}\; \gamma(A)\geq k\}.$$

We state now the symmetric mountain pass lemma of Kajikia (see Theorem 1 in \cite{K}).
\begin{teo}\label{t3}
Assume $X$ is an infinite dimensional Banach space and $\Lambda\in C^1(X,\RR)$ satisfies conditions
(A1) and (A2) below.
\smallskip

\noindent (A1) $\Lambda(u)$ is even, bounded from below, $\Lambda(0)=0$ and $\Lambda(u)$ satisfies
the Palais-Smale condition (i.e., any sequence $\{u_n\}$ in $X$ such that $\{\Lambda(u_n)\}$ is bounded
and $\Lambda^{'}(u_n)\rightarrow 0$ in $X^\star$ as $n\rightarrow\infty$ has a convergent subsequence);
\smallskip

\noindent (A2) For each $k\in\NN$, there exists an $A_k\in\Gamma_k$ such that $\sup_{u\in A_k}\Lambda(u)<0$.
\smallskip

Under the above assumptions, either (i) or (ii) below holds true.
\smallskip

\noindent (i) There exists a sequence $\{u_n\}$ such that $\Lambda^{'}(u_n)=0$, $\Lambda(u_n)<0$ and
$\{u_n\}$ converges to zero;
\smallskip

\noindent (ii) There exist two sequences $\{u_n\}$ and $\{v_n\}$ such that $\Lambda^{'}(u_n)=0$,
$\Lambda(u_n)=0$, $u_n\neq 0$, $\lim_{n\rightarrow\infty}u_n=0$, $\Lambda^{'}(v_n)=0$,
$\Lambda(v_n)=0$, and $v_n$ converges to a non-zero limit.
\end{teo}

In order to apply Theorem \ref{t3} to the functional $J_\lambda$ we prove two auxiliary results.

\begin{lemma}\label{l3}
The functional $J_\lambda$ satisfies condition (A1) from Theorem \ref{t3}.
\end{lemma}
\proof
Clearly, $J_\lambda(u)=J_\lambda(-u)$ for any $u\in E$, i.e. $J_\lambda$ is even, and $J_\lambda(0)=0$.
On the other hand, since by relation \eq{7} we have
$$\int_\Omega\Phi(|\nabla u(x)|)\;dx\geq\|u\|^{p^0},\;\;\;\forall\;
u\in E\;{\rm with}\;\|u\|<1,$$
while by Lemma C.9 in \cite{Clem2} we have
$$\int_\Omega\Phi(|\nabla u(x)|)\;dx\geq\|u\|^{p_0},\;\;\;\forall\;
u\in E\;{\rm with}\;\|u\|>1,$$
we deduce that
\begin{equation}\label{stea3}
\int_\Omega\Phi(|\nabla u(x)|)\;dx\geq\alpha(\|u\|),\;\;\;\forall\;
u\in E,
\end{equation}
where $\alpha:[0,\infty)\rightarrow\RR$, $\alpha(t)=t^{p^0}$ if $t<0$ and $\alpha(t)=t^{p_0}$ if $t>1$.

By Remark 1, the space $E$ is continuously embedded in $L^{q^{\pm}}(\Omega)$. Thus, there exist two positive 
constants
$d_1$ and $d_2$ such that
\begin{equation}\label{stea4}
\int_\Omega|u|^{q^+}\;dx\leq d_1\|u\|^{q^+},\;\;\;\int_\Omega|u|^{q^-}\;dx\leq d_2\|u\|^{q^-},\;\;\;\forall\;u\in E.
\end{equation}
Combining relations \eq{stea3} and \eq{stea4} we get
$$J_\lambda(u)\geq\alpha(\|u\|)-\frac{d_1\lambda}{q^-}\|u\|^{q^+}-\frac{d_2\lambda}{q^-}\|u\|^{q^-},\;\;\;\forall\;u\%
in E.$$
Since by relation \eq{stea2} we have $q^+<p_0$ the above relation
shows that $J_\lambda$ is bounded from below.

Next, we show that $J_\lambda$ satisfies the Palais-Smale condition. Let $\{u_n\}$ be a
sequence in $E$ such that $\{J_\lambda(u_n)\}$ is bounded and $J^{'}(u_n)\rightarrow 0$ in $E^{'}$, as 
$n\rightarrow\infty$.
We show that $\{u_n\}$ is bounded in $E$. Assume by
contradiction the contrary. Then, passing eventually to a subsequence,
still denoted by $\{u_n\}$, we may assume that $\|u_n\|
\rightarrow\infty$ as $n\rightarrow\infty$. Thus we may consider that
$\|u_n\|>1$ for any integer $n$.

By our assumptions, there is a positive constant $M$ such that for all $n$ large enough we have
\begin{eqnarray*}
M+1+\|u_n\|&\geq&J_\lambda(u_n)-\frac{1}{q^-}\langle J^{'}(u_n),u_n\rangle\\
&=&\int_{\Omega}\Phi(|\nabla u_n|)\;dx-\lambda
\int_\Omega\frac{1}{q(x)}|u_n|^{q(x)}\;dx-\frac{1}{q^-}\cdot
\int_\Omega\phi(|\nabla u_n(x)|)|\nabla u_n(x)|\;dx+\\
&&\frac{\lambda}{q^-}
\int_\Omega|u_n|^{q(x)}\;dx\\
&\geq&\int_{\Omega}\Phi(|\nabla u_n|)\;dx-\frac{1}{q^-}\cdot
\int_\Omega\phi(|\nabla u_n(x)|)|\nabla u_n(x)|\;dx\\
&\geq&\left(1-\frac{p^0}{q^-}\right)\int_{\Omega}\Phi(|\nabla u_n|)\;dx\\
&\geq&\left(1-\frac{p^0}{q^-}\right)\|u_n\|^{p_0}.
\end{eqnarray*}
Since $p_0>1$, letting
$n\rightarrow\infty$ we obtain a contradiction. It follows that
$\{u_n\}$ is bounded in $E$.
Similar arguments as those used in the end of the proof of Theorem \ref{t1} imply that, up
to a subsequence, $\{u_n\}$ converges
strongly in $E$.

The proof of Lemma \ref{l3} is complete.
\endproof

\begin{lemma}\label{l4}
The functional $J_\lambda$ satisfies condition (A2) from Theorem \ref{t3}.
\end{lemma}
\proof
We construct a sequence of subsets $A_k\in\Gamma_k$ such that $\sup_{u\in A_k}J_\lambda(u)<0$,
for each $k\in\NN$.

Let $x_1\in\Omega$ and $r_1>0$ be such that $\overline{B_{r_1}(x_1)}\subset\Omega$ and
$|\overline{B_{r_1}(x_1)}|<|\Omega|/2$. Consider $\theta_1\in C_0^\infty(\Omega)$ be a function
with ${\rm supp}(\theta_1)=\overline{B_{r_1}(x_1)}$.

Define $\Omega_1=\Omega\setminus\overline{B_{r_1}(x_1)}$.

Next, let $x_2\in\Omega$ and $r_2>0$ be such that $\overline{B_{r_2}(x_2)}\subset\Omega_1$ and
$|\overline{B_{r_2}(x_2)}|<|\Omega_1|/2$. Consider $\theta_2\in C_0^\infty(\Omega)$ be a function
with ${\rm supp}(\theta_2)=\overline{B_{r_2}(x_2)}$.

Continuing the process described above we can construct by recurrence a sequence of functions
$\theta_1$, $\theta_2$,..., $\theta_k\in C_0^\infty(\Omega)$ such that ${\rm supp}(\theta_i)\neq
{\rm supp}(\theta_j)$ if $i\neq j$ and $|{\rm supp}(\theta_i)|>0$ for any $i,\;j\in\{1,...,k\}$.

We define the finite dimensional subspace of $E$,
$$F={\rm span}\{\theta_1,\theta_2,...,\theta_k\}.$$
Clearly, ${\rm dim}F=k$ and $\int_\Omega|\theta|^{q(x)}\;dx>0$, for any $\theta\in F\setminus\{0\}$.
We denote by $S_1$ the unit sphere in $E$, i.e. $S_1=\{u\in E;\;\|u\|=1\}$. For any number $t\in(0,1)$
we define the set
$$A_k(t)=t\cdot(S_1\cap F).$$
Since for any bounded symmetric neighborhood $\omega$ of the origin in $\RR^k$ there holds $\gamma(\partial\omega)=k$
(see Proposition 5.2 in \cite{S}) we deduce that $\gamma(A_k(t))=k$ for any $t\in(0,1)$.

Finally, we show that for each integer $k$ there exists $t_k\in(0,1)$ such that
$$\sup_{u\in A_k(t_k)}J_\lambda(u)<0.$$
For any $t\in(0,1)$ we have
\begin{eqnarray*}
\sup_{u\in A_k(t)}J_\lambda(u)&\leq&\sup_{\theta\in S_1\cap F}J_\lambda(t\theta)\\
&=&\sup_{\theta\in S_1\cap F}\left\{\int_{\Omega}\Phi(t|\nabla\theta|)\;dx-\lambda
\int_\Omega\frac{1}{q(x)}t^{q(x)}|\theta|^{q(x)}\;dx\right\}\\
&\leq&\sup_{\theta\in S_1\cap F}\left\{t^{p_0}\int_{\Omega}\Phi(|\nabla\theta|)\;dx-\frac{\lambda t^{q^+}}{q^+}
\int_\Omega|\theta|^{q(x)}\;dx\right\}\\
&=&\sup_{\theta\in S_1\cap F}\left\{t^{p_0}\left(1-\frac{\lambda}{q^+}\cdot\frac{1}{t^{p_0-q^+}}\cdot
\int_\Omega|\theta|^{q(x)}\;dx\right)\right\}
\end{eqnarray*}
Since $S_1\cap F$ is compact we have $m=\min_{\theta\in S_1\cap F}\int_\Omega|\theta|^{q(x)}\;dx>0$. Combining
that fact with the information given by relation \eq{stea2}, that is $p_0>q^+$, we deduce that we can choose
$t_k\in(0,1)$ small enough such that
$$1-\frac{\lambda}{q^+}\cdot\frac{1}{t^{p_0-q^+}}\cdot m<0.$$
The above relations yield
$$\sup_{u\in A_k(t_k)}J_\lambda(u)<0.$$
The proof of Lemma \ref{l4} is complete.
\endproof

{\sc Proof of Theorem \ref{t2}.} Using Lemmas \ref{l3} and \ref{l4} we deduce that we can apply Theorem \ref{t3}
to the functional $J_\lambda$. So, there exists a sequence $\{u_n\}\subset E$ such that $J^{'}(u_n)=0$,
for each $n$, $J_\lambda(u_n)\leq 0$ and $\{u_n\}$ converges to zero in $E$.

The proof of Theorem \ref{t2} is complete.  \endproof
\section{Examples}
In this section we point out two concrete examples of problems to which we can apply the main results of
this paper.

\noindent{\sc Example 1.} We consider the problem
\begin{equation}\label{ex1}
\left\{\begin{array}{lll}
-{\rm div}(\log(1+|\nabla u|^r)|\nabla u|^{p-2}\nabla
u)=\lambda|u|^{q(x)-2}u,
&\mbox{for}& x\in\Omega\\
u=0, &\mbox{for}& x\in\partial\Omega,
\end{array}\right.
\end{equation}
where $p$ and $r$ are real numbers such that $1<p$, $r$, $N>p+r$ and $q(x)$ is a continuous function on 
$\overline\Omega$
such that $1<q(x)$ for all $x\in\overline\Omega$ and furthermore
$$\inf_\Omega q(x)<p\;\;\;{\rm and}\;\;\;\sup_\Omega q(x)<\frac{Np}{N-p}.$$
In this case we have
$$\phi(t)=\log(1+|t|^r)\cdot|t|^{p-2}t,\qquad \mbox{for all}\ t\in\RR$$
and
$$\Phi(t)=\int_0^t\phi(s),\qquad \mbox{for all}\ t\in\RR.$$
Clearly, $\phi$ is an odd, increasing
homeomorphism of $\RR$ into $\RR$, while $\Phi$ is convex and even
on $\RR$ and increasing from $\RR_+$ to $\RR_+$.

By Example 2 on p. 243 in \cite{Clem2} we know that
$$p_0=p\;\;\;{\rm and}\;\;\;p^0=p+r$$
and thus relation \eq{6} in Theorem \ref{t1} is satisfied. On the other hand, by Proposition 1 in \cite{JMAA}
(see also \cite{CRAS}) we deduce that relations \eq{2} and \eq{6biss} are fulfilled. Thus, we verified that
we can apply Theorem \ref{t1} in order to find out that there exists $\lambda^\star>0$ such that
any $\lambda\in(0,\lambda^\star)$ is an eigenvalue for problem \eq{ex1}.

\noindent{\sc Example 2.} We consider the problem
\begin{equation}\label{ex2}
\left\{\begin{array}{lll}
-{\rm div}\left(\di\frac{|\nabla u|^{p-2}\nabla u}{\log(1+|\nabla u|)}\right)=\lambda|u|^{q(x)-2}u, &\mbox{for}&
x\in\Omega\\
u=0, &\mbox{for}& x\in\partial\Omega\,,
\end{array}\right.
\end{equation}
where $p$ is a real number such that $p>N+1$ and $q\in C(\overline\Omega)$ satisfies
$1<q(x)<p-1$ for any $x\in\overline\Omega$. In this case we have
$$\phi(t)=\frac{|t|^{p-2}}{\log(1+|t|)}t$$
and
$$\Phi(t)=\int_0^t\phi(s)\;ds,$$
is an increasing continuous function from $\RR^+$ to $\RR^+$, with $\Phi(0)=0$ and such that the function
$\Phi(\sqrt{t})$ is convex. By Example 3 on p. 243 in \cite{Clem2} we have
$$p_0=p-1<p^0=p=\liminf_{t\rightarrow\infty}\frac{\log(\Phi(t))}{\log(t)}.$$
Thus, conditions \eq{2}, \eq{stea2} and \eq{66} from Theorem
\ref{t2} and Remark 2 are verified. We deduce that every
$\lambda>0$ is an eigenvalue for problem \eq{ex2}. Moreover, for
each $\lambda>0$ there exists a sequence of eigenvectors $\{u_n\}$
such that $\lim_{n\rightarrow\infty}u_n=0$ in $W_0^1L_\Phi(\Omega)$.

\medskip
{\bf Acknowledgements.}
The authors have been
supported by Grant CNCSIS PNII--79/2007 {\it ``Procese
Neliniare Degenerate \c si Singulare"}.

\end{document}